\documentclass[11pt,bezier]{article}
\usepackage{amsmath}
\usepackage{amsfonts,amsthm,amssymb}
\usepackage{amsfonts}
\usepackage{tikz, color}
\usepackage{graphics,epstopdf}
\usepackage{setspace}
\usepackage{enumerate}
\usepackage{graphicx}

\usepackage{graphics}
\usepackage{amssymb}

\newtheorem{lem}{Lemma}[section]
\newtheorem{thm}[lem]{Theorem}

\newtheorem{conj}{Conjecture}

\theoremstyle{definition}

\begin{document}
\title{Edge-connectivity keeping trees in $k$-edge-connected graphs\footnote{The research is supported by National Natural Science Foundation of China (11861066).}}
\author{Qing Yang, Yingzhi Tian\footnote{Corresponding author. E-mail: tianyzhxj@163.com (Y. Tian), qingyangmath@sina.com (Q. Yang).} \\
{\small College of Mathematics and System Sciences, Xinjiang
University, Urumqi, Xinjiang, 830046, PR China}}

\date{}

\maketitle

\noindent{\bf Abstract } Mader [J. Combin. Theory Ser. B 40 (1986) 152-158] proved that every $k$-edge-connected graph $G$ with minimum degree at least $k+1$ contains a vertex $u$ such that $G-\{u\}$ is still $k$-edge-connected. In this paper, we prove that every $k$-edge-connected graph $G$ with minimum degree at least $k+2$ contains an edge $uv$ such that $G-\{u,v\}$ is $k$-edge-connected for any positive integer $k$. In addition, we show that for any tree $T$ of order $m$, every $k$-edge-connected graph $G$ with minimum degree greater than $4(k+m)^2$ contains a subtree $T'$ isomorphic to $T$ such that $G-V(T')$ is $k$-edge-connected.


\section{Introduction}

Throughout this paper, graph always means a finite, undirected graph without multiple edges and loops. Let $G$ be a graph with \emph{vertex set} $V(G)$, \emph{edge set} $E(G)$ and \emph{order} $|V(G)|$. For a vertex $v\in V(G)$, the \emph{neighborhood} of $v$ is $N_{G}(v)$. For a vertex set $S\subseteq V(G)$, the \emph{neighborhood} $N_{G}(S)$  of $S$ in $G$ is $(\cup_{v\in S}N_{G}(v))\setminus S$. The \emph{degree} $d_G(v)$ of $v$ in $G$ is $|N_{G}(v)|$. The \emph{minimum degree} $\delta(G)$ of $G$ is min$_{v\in V(G)}d_G(v)$. The \emph{induced subgraph} $G[S]$, is the graph with vertex set $S$, where two vertices in $S$ are adjacent if and only if they are adjacent in $G$. And $G-S$ is the induced graph $G[V(G)\backslash S)]$. For an edge set $U\subseteq E(G)$, let $G-U$ be the subgraph of $G$ with vertex set $V(G)$ and edge set $E(G)\backslash U$. For graph-theoretical terminologies and notation not defined here, we follow \cite{Bondy}.

The \emph{connectivity} of $G$, denoted by $\kappa(G)$, is the minimum size of a vertex set $S\subseteq V(G)$ such that $G-S$ is disconnected or has only one vertex. The \emph{edge-connectivity} of $G$, denoted by $\kappa'(G)$, is the minimum size of an edge set $U\subseteq E(G)$ such that $G-U$ is disconnected. The graph $G$ is said to be \emph{$k$-connected} (respectively, \emph{$k$-edge-connected}) if $\kappa(G)\geq k$ (respectively, $\kappa'(G)\geq k$). Note that a trivial graph is both 1-connected and 1-edge-connected, but is not $k$-connected and $k$-edge-connected for any $k > 1$.

In 1972, Chartrand, Kaigars and Lick \cite{Chartrand} proved that there is a redundant vertex in a $k$-connected graph with minimum degree at least $\lfloor\frac{3k}{2}\rfloor$.

\begin{thm} (Chartrand et al.\cite{Chartrand}) Every $k$-connected graph $G$ with $\delta(G)\geq \lfloor\frac{3k}{2}\rfloor$ contains a vertex $u$ such that $\kappa (G-\{u\})\geq k$.
\end{thm}

In 2008, Fujita and Kawarabayashi \cite{Fujita} showed that every $k$-connected graph $G$ with $\delta(G)\geq \lfloor\frac{3k}{2}\rfloor+2$ contains an edge $uv$ such that $\kappa (G-\{u,v\})\geq k$. Furthermore, they conjectured that every $k$-connected graph with $\delta(G)\geq \lfloor\frac{3k}{2}\rfloor+f_{k}(m)-1$ contains a connected subgraph $H$ of order $m$ such that $\kappa (G-V(H))\geq k$, where $f_{k}(m)$ is nonnegative. Mader \cite{Mader1} confirmed the conjecture and proved that $f_{k}(m)=m$. In addition, the connected subgraph $H$ can be chosen as a path. Meanwhile, Mader conjectured that the result would hold even if we replace the path by any tree with the same order.

\begin{conj} (Mader \cite{Mader1}) For any tree $T$ with order $m$, every $k$-connected graph $G$ with $\delta(G)\geq \lfloor\frac{3k}{2}\rfloor+m-1$ contains a tree $T'\cong T$ such that $\kappa (G-V(T'))\geq k$.
\end{conj}

Actually, the theorem in \cite{Diwan} implies Conjecture 1 holds for $k=1$. When $k=2$, Conjecture 1 holds for star, path-star, spider and caterpillar et al., see [5-7, 9, 14-15] for references. In 2022, Hong and Liu \cite{Hong} confirmed Conjecture 1 for $k=2$ and $k=3$. Conjecture 1 remains open for $k\geq4$. However, Mader \cite{Mader2} showed that Conjecture 1 holds if $\delta(G)\geq2(k-1+m)^2+m-1$.

In this paper, we consider a similar problem in the $k$-edge-connected graphs. Mader \cite{Mader0} showed that there is a redundant vertex in a $k$-edge-connected graph with minimum degree at least $k+1$.

\begin{thm} (Mader \cite{Mader0}) Every $k$-edge-connected graph $G$ with $\delta(G)\geq k+1$ contains a vertex $u$ such that $\kappa'(G-\{u\})\geq k$.
\end{thm}

In this paper, we prove that there exists a pair of redundant adjacent vertices in a $k$-edge-connected graph $G$ with $\delta(G)\geq k+2$ for any positive integer $k$. Meanwhile, we also obtain that for any tree $T$ of order $m$, every $k$-edge-connected graph $G$ with $\delta(G)>4(k+m)^2$ contains a subtree $T'\cong T$ such that $\kappa'(G-V(T'))\geq k$. In the next section, we will introduce some definitions and lemmas that will be used to prove the main results. In Section 3, we present the main results. Concluding remarks will be given in the last section.

\section{Preliminaries}

Let $U$ be an edge-cut of $G$ and $F$ be the union of at least one component, but not all components of $G-U$. Then $F$ is called a \emph{semifragment} of $G$ to $U$. And $F'=G-V(F)$ is called a \emph{complementary semifragment} of $F$ in $G$ to $U$. In fact, $F'$ is also a semifragment of $G$ to $U$. In addition, if $U$ is a minimum edge-cut of $G$, then $G-U$ has exactly two components $F$ and $F'$. Each of $F$ and $F'$ is a \emph{fragment} of $G$ to $U$, and $F'$ is a complementary fragment of $F$ in $G$ to $U$. Sometimes, we simply say $F$ is a fragment of $G$.

For two vertex set $V_1, V_2\subseteq G$, the edge set between $V_1$ and $V_2$ in $G$ is denoted by $E_{G}[V_1,V_2]$ and the number of all edges between $V_1$ and $V_2$ in $G$ is denoted by $d_{G}(V_{1}, V_{2})$, which is $|E_{G}[V_1,V_2]|$ or $\sum_{u\in V_1 }|N_{G}(u)\cap V_2|$. For convenience, this paper may uses $G[V_1]$ and $G[V_2]$ for $V_1$ and $V_2$, respectively. For examples, we use $G[V_1]\cap G[V_2]$, $G[V_1]\cup G[V_2]$ and $d_{G}(G[V_1], G[V_2])$ for $V_1\cap V_2$, $V_1\cup V_2$ and $d_{G}(V_1, V_2)$, respectively.

\begin{lem} Let $e$ and $e_{1}$ be two nonadjacent edges of $G$, and, let $U$ and $U_1$ be minimum edge-cuts of $G-V(e)$ and $G-V(e_{1})$, respectively. Assume that $F$ is a fragment of $G-V(e)$ to $U$ and $F_{1}$ is a fragment of $G-V(e_{1})$ to $U_{1}$. And we denote $F'=G-V(e)-V(F)$ and $F_{1}'=G-V(e_{1})-V(F_{1})$. Furthermore, we assume $e\in E(F_{1})$, $e_{1}\in E(F')$ and $F\cap F_{1}\neq\emptyset$ (see Figure 1). Then we have the following properties.

($\alpha$) If $F'\cap F'_{1}\neq\emptyset$ holds, then $F\cap F_{1}$ is a fragment of $G-V(e)$.

($\beta$) If $F'\cap F'_{1}=\emptyset$ holds, then $|V(F')|<|V(F_{1})|$.
 \end{lem}
\begin{figure}[htbp]
\centering
 \includegraphics[width=0.25\textwidth]{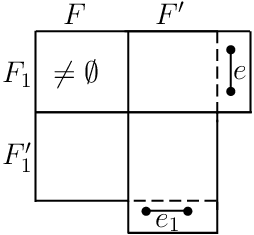}\\
\caption{Illustration for Lemma 2.1}
\end{figure}

\noindent{\bf Proof.} ($\alpha$) Since $F\cap F_{1}\neq\emptyset$, we have
\begin{equation}
\begin{aligned}
d_{G-V(e)}(F\cap F_{1},F'\cup F'_{1})&\geq \kappa'(G-V(e))\\
&=d_{G}(F, F') \\
&= d_{G}(F\cap F_{1}, F')+d_{G}(F\cap F'_{1}, F').
\end{aligned}
\end{equation}
Thus
\begin{equation}
\begin{aligned}
d_{G}(F \cap F_{1}, F\cap F'_{1})&\geq d_{G}(F\cap F'_{1} , F')
\end{aligned}
\end{equation}
by $d_{G-V(e)}(F\cap F_{1},F'\cup F'_{1})=d_{G}(F\cap F_{1}, F')+d_{G}(F\cap F_{1}, F\cap F'_{1} ).$

If $F'\cap F'_{1}\neq\emptyset$ holds, then
\begin{equation}
\begin{aligned}
d_{G-V(e_1)}(F'\cap F'_{1},F\cup F_{1})&\geq \kappa'(G-V(e_1))\\
&=d_{G}(F_1, F_1') \\
&= d_{G}(F'\cap F'_{1}, F_1)+d_{G}(F\cap F'_{1}, F_1).
\end{aligned}
\end{equation}
Thus
\begin{equation}
\begin{aligned}
d_{G}(F\cap F'_{1}, F'\cap F'_{1})&\geq d_{G}(F\cap F'_{1}, F_{1})
\end{aligned}
\end{equation}
by $d_{G-V(e_1)}(F'\cap F'_{1},F\cup F_{1})=d_{G}(F'\cap F'_{1}, F_1)+d_{G}(F'\cap F'_{1}, F\cap F'_{1}).$
By inequalities (2) and (4), we get
\begin{equation}
\begin{aligned}
d_{G}(F\cap F_{1} , F\cap F'_{1} )&\geq d_{G}(F\cap F'_{1}, F')\\
&\geq d_{G}( F\cap F'_{1}, F'\cap F'_{1})\\
&\geq d_{G}(F\cap F'_{1}, F_{1})\\
&\geq d_{G}(F\cap F_{1} , F\cap F'_{1} ),
\end{aligned}
\end{equation}
that is, $d_{G}(F \cap F_{1}, F\cap F'_{1})=d_{G}(F\cap F'_{1} , F')$. Moreover, we have
\begin{equation}
\begin{aligned}
d_{G-V(e)}(F\cap F_{1},G-V(e)-(F\cap F_{1}))
&=d_{G-V(e)}(F\cap F_{1},F'\cup F'_{1})\\
&=d_{G}(F\cap F_{1}, F')+d_{G}(F\cap F_{1}, F\cap F'_{1})\\
&= d_{G}( F\cap F_{1}, F')+d_{G}(F\cap F'_{1} , F')\\
&= d_{G}(F, F')\\
&= \kappa'(G-V(e)).
\end{aligned}
\end{equation}
Thus, $U'=E_{G-V(e)}[F\cap F_{1},G-V(e)-(F\cap F_{1})]$ is a minimum edge-cut of $G-V(e)$ and $F\cap F_{1}$ is a fragment of $G-V(e)$ to $U'$.

($\beta$) If $F'\cap F'_{1}=\emptyset$ holds, then $F'-V(e_{1})\subseteq F_1$. Since $V(e)\in V(F_1)$, $V(e)\notin V(F')$ and $F\cap F_{1} \neq\emptyset$, we have $|V(F')|<|V(F_{1})|$. $\Box$

\begin{lem} Let $G$ be a $k$-edge-connected graph with $\delta(G)\geq k+2$ and let $e_{0}$ be an edge of graph $G$ such that $\kappa'(G-V(e_{0}))<k$. Assume that $U_0$ is a minimum edge-cut of $G-V(e_0)$ and $F_0$ is a fragment of $G-V(e_0)$ to $U_0$. Then there exist an edge $e_1\in E(G[V(F_0)\cup V(e_0)])$ with $\kappa'(G-V(e_{1}))<k$ and a fragment $F_1$ of $G-V(e_{1})$ with $F_1\subseteq G[V(F_0)\cup V(e_0)]$ such that the following property holds.\begin{quote}
For any $e\in E(F_1)$, if $\kappa'(G-V(e))< k$ holds, then any fragment $F$ of $G-V(e)$ with $e_1\notin E(F)$ satisfies $F\cap F_{1}=\emptyset$. See Figure 2 for an illustration.
\end{quote}
 \end{lem}
\begin{figure}[htbp]
\centering
 \includegraphics[width=0.26\textwidth]{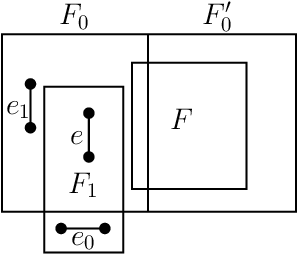}\\
\caption{Illustration for Lemma 2.2}
\end{figure}
\noindent{\bf Proof.} Let $E_0=\{e\in E(G[V(F_0)\cup V(e_0)]): \kappa'(G-V(e))<k\}$ and for any $e\in E_{0}$, let
\begin{center}
$\mathcal{U}_e=\{U$: $U$~is~a~minimum~edge-cut~of~$G-V(e)$\}.
\end{center}
For any $U_e\in \mathcal{U}_e$, we denote
\begin{center}
$\mathcal{F}_{U_e}=\{F$: $F$~is~a~fragment~of~$G-V(e)$~to~$U_e$~with~$F\subseteq G[V(F_0)\cup V(e_0)]$\}.
\end{center}
Furthermore, let $\mathcal{F}=\cup_{e\in E_0}(\cup_{U_e\in \mathcal{U}_e}\mathcal{F}_{U_e})$. Then $e_0\in E_0$, $U_0\in \mathcal{U}_{e_0}$ and $F_0\in \mathcal{F}_{\mathcal{U}_{e_0}}$. Thus $\mathcal{F}\neq\emptyset$. Pick $F_1\in\mathcal{F}$ so that $|V(F_1)|$ is as small as possible. Assume that $F_1$ is a fragment of $G-V(e_1)$ to $U_{e_1}$ and $F'_1=G-V(e_1)-V(F_1)$ is the complementary fragment of $F_1$ in $G-V(e_1)$ to $U_{e_1}$, where $e_1\in E_0$, $U_{e_1}\in \mathcal{U}_{e_1}$ and $F_1\in \mathcal{F}_{\mathcal{U}_{e_1}}$. 
Note that  $\delta(F_1)\geq k+2-(k-1)-2=1$. Assume $e\in E(F_1)$ satisfies $$\kappa'(G-V(e))<k.$$ Then $e\in E_0$ by $F_1\subseteq G[V(F_0)\cup V(e_0)]$. Let $U$ be a minimum edge-cut of $G-V(e)$, $F$ be a fragment of $G-V(e)$ to $U$ with $e_1\notin E(F)$ and $F'=G-V(e)-V(F)$ be the complementary fragment of $F$ in $G-V(e)$ to $U$. Then  $e_1\in E(F')$.

Suppose $F\cap F_{1}\neq\emptyset$. By Lemma 2.1 ($\alpha$), if $F'\cap F'_{1}\neq\emptyset$ hods, then $F\cap F_{1}$ is a fragment of $G-V(e)$. Since $F_1\subseteq G[V(F_0)\cup V(e_0)]$, $\kappa'(G-V(e))<k$ and $e\in E(F_1)$, we get $F\cap F_{1} \in \mathcal{F}$. By $e\in E(F_1)$ and $e\notin E(F)$, we have $|F\cap F_{1}|< |V(F_{1})|$, which contradicts the choice of $F_1$. By Lemma 2.1 ($\beta$), if $F'\cap F'_{1}=\emptyset$ holds, then $|V(F')|<|V(F_{1})|$. Since $ F'\cap F'_{1}=\emptyset$ and $e_1\in E(F')$, we obtain $$F'\subseteq G[V(F_{1})\cup V(e_1)]\subseteq G[V(F_{0})\cup V(e_0)].$$  Thus $F'\in \mathcal{F}$, also a contradiction. $\Box$

\section{Main results}

\begin{thm} For any positive integer $k$, every graph $G$ with $\kappa'(G)\geq k$ and $\delta(G)\geq k+2$ contains an edge $e$ such that $\kappa'(G-V(e))\geq k$.
\end{thm}

\noindent{\bf Proof.} By contrary, assume that $\kappa'(G-V(e))< k$ for any $e\in E(G)$. Then we have the following claims.

\noindent{\bf Claim 1.} There exist an edge $e_1\in E(G)$ and a fragment $F_1$ of $G-V(e_1)$ such that $$2\leq|V(F_1)|<k.$$

%$U_1$ is a minimum edge of $G-V(e_1)$ and $F_1$ is a fragment of $G-V(e_1)$ to $U_1$. By $\delta(F_1)\geq k+2-(k-1)-2=1$, we get $|V(F_1)|\geq 2$. Since $\kappa'(G-V(e))< k$ for any $e\in E(F_1)$, by Lemma 2.2, there exists a minimum edge-cut $U$ such that $F\cap F_{1}=\emptyset$ for any fragment $F$ of $G-V(e)$ to $U$ with $e_1\in E(F')$, where $F'=G-V(e)-V(F)$.

Let $e_0\in E(G)$. Assume that $U_0$ is a minimum edge-cut of $G-V(e_0)$ and $F_0$ is a fragment of $G-V(e_0)$ to $U_0$. Since $\kappa'(G-V(e_0))< k$, by Lemma 2.2, there exist $e_1\in E(G[V(F_0)\cup V(e_0)])$ and a fragment $F_1\subseteq G[V(F_0)\cup V(e_0)]$ of $G-V(e_1)$ such that for any $e\in E(F_1)$ and any fragment $F$ of $G-V(e)$ with $e_1\notin E(F)$, we have $F\cap F_{1}=\emptyset$. Let $F'=G-V(e)-V(F)$ and $F_{1}'=G-V(e_{1})-V(F_{1})$. Thus $F\subseteq F'_1$. Since $\kappa'(G-V(e))< k$, we have $d_{G}(V(e),F)>0$. By $F\subseteq F'_1$, we get $d_{G}(V(e),F'_{1})>0$. If there exists a vertex $u\in V(F_1)$ such that $N_{G}(u)\cap V(F'_1)=\emptyset$, then by $d_{G}(V(e),F'_1)>0$ for any $e\in E(F_1)$, we have $N_{G}(v)\cap V(F'_1)\neq \emptyset$ for any $v\in N_{G}(u)\cap V(F_1)$. Since
\begin{equation}
\begin{aligned}
|N_{G}(u)\cap V(F_1)|&=d_{G}(u)-|N_{G}(u)\cap V(F'_1)|-|N_{G}(u)\cap V(e_1)|\\
&\geq k+2-0-2\\
&=k,
\end{aligned}
\end{equation}
we have $k\leq |N_{G}(u)\cap V(F_1)| \leq d_{G}(F_1,F'_1)=\kappa'(G-V(e_1))<k$, a contradiction. Thus, for any $w\in V(F_1)$, we obtain $N_{G}(w)\cap V(F'_1)\neq\emptyset$. Since $d_{G}(F_1,F'_1)<k$, we have $|V(F_1)|<k$. Thus Claim 1 holds.

By claim 1, the following claim is established.
%By Claim 1, we have the following claim.

\noindent{\bf Claim 2.} $k\geq3.$

%By Claim 1, we have $2\leq|V(F_1)|<k$, that is, $k\geq3.$

%\begin{figure}[htbp]
%\centering
%  % Requires \usepackage{graphicx}
% \includegraphics[width=0.25\textwidth]{23-11-22-3.eps}\\
%\caption{Partition of $G$}
%\end{figure}

Since $G$ is a simple graph with $\delta(G)\geq k+2$, we have $$d_{G}(F_1, F'_1)\geq|V(F_1)|(k-|V(F_1)|+1)$$ by \begin{equation}
\begin{aligned}
|N_{G}(z)\cap V(F'_1)|&=d_{G}(z)-|N_{G}(z)\cap V(F_1)|-|N_{G}(z)\cap V(e_1)|\\
&\geq k+2-(|V(F_1)|-1)-2\\
&=k-|V(F_1)|+1
\end{aligned}
\end{equation} 
for any $z\in V(F_1)$.  Thus $d_{G}(F_1, F'_1)\geq|V(F_1)|(k-|V(F_1)|+1)\geq 2(k-1)$ by $2\leq|V(F_1)|<k$. 
This, together with $d_{G}(F_1, F'_1)< k$, we obtain $k<2$, which contradicts to $k\geq3$ (by Claim 2). $\Box$

The following lemma will be used to prove that there exists a redundant subtree in $k$-edge-connected graphs with minimum degree greater than $4(k+m)^2$.

\begin{lem} (Thomassen \cite{Thomassen}) Any graph $G$ with $\delta(G)> 4k^2$ contains a $k$-connected subgraph $H$ with $|V(H)|>4k^2$ such that $|N_G(G-V(H))|\leq 2k^2$.
\end{lem}

\begin{thm} For any tree $T$ with order $m$, every $k$-edge-connected graph $G$ with $\delta(G)>4(k+m)^2$ contains a tree $T'\cong T$ such that $G-V(T')$ is still $k$-edge-connected.
\end{thm}
\noindent{\bf Proof.} By Lemma 3.2, there exists a $(k+m)$-connected subgraph $H$ with $|V(H)|>4(k+m)^2$ in $G$ such that $|N_G(G-V(H))|\leq 2(k+m)^2$. Let $\bar{H}=G-V(H)$. Since $\delta(G)>4(k+m)^2$, we obtain $\delta(H-N_G(\bar{H}))\geq m$, thus there is a subtree $T'\cong T$ in $H-N_G(\bar{H})$. If $G-V(T')$ is $k$-edge-connected, then the result holds. So, we assume $\kappa'(G-V(T'))\leq k-1$.

Let $U$ be a minimum edge-cut of $G-V(T')$. Then $|U|\leq k-1$. Let $F$ be a fragment of $G-V(T')$ to $U$ and $F'=G-V(T')-V(F)$. Note that $F'$ is also a fragment of $G-V(T')$ to $U$. We give some notations as follows (see Figure 3 for an illustration).
$$H_1=F\cap H,~H_2=F'\cap H;$$
$$ \bar{H}_1=F\cap \bar{H},~ \bar{H}_2 =F'\cap \bar{H};$$
$$D_1=V(U)\cap F,~D_2=V(U)\cap F'.$$

\begin{figure}[htbp]
\centering
 \includegraphics[width=0.35\textwidth]{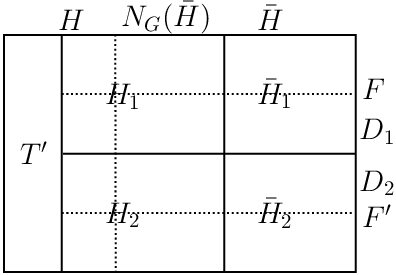}\\
\caption{Illustration for the proof of Theorem 3.3}
\end{figure}

By $|U|\leq k-1$, we get $|D_1|\leq k-1$ and $|D_2|\leq k-1$. Since $H$ is $(k+m)$-connected, $H-D_1\cup V(T')$ is connected, thus \begin{equation}
\begin{aligned}
H_1\backslash D_1=\emptyset~ or ~H_2=\emptyset.
\end{aligned}
\end{equation} 
And by $H-D_2\cup V(T')$ is connected, we have
\begin{equation}
\begin{aligned}
H_2\backslash D_2=\emptyset~ or ~H_1=\emptyset.
\end{aligned}
\end{equation} 
If both $H_1\backslash D_1=\emptyset$ and $H_2\backslash D_2=\emptyset$, then $|V(H)\backslash V(T')|\leq|D_1\cup D_2|\leq 2(k-1)$, which is a contradiction to $|V(H)|>4(k+m)^2$. Thus \begin{equation}
\begin{aligned}
H_1\backslash D_1\neq\emptyset~ or ~H_2\backslash D_2\neq\emptyset.
\end{aligned}
\end{equation}

If $H_1\neq\emptyset$ holds, then $H_2\backslash D_2=\emptyset$ by (10). Suppose that $H_2\neq\emptyset$ holds. Then $H_1\backslash D_1=\emptyset$ by (9), a contradiction to (11). Thus $H_2=\emptyset$. Since $F'$ is a fragment of $G-V(T')$ to $U$, we have $F'\neq\emptyset$, thus $\bar{H}_2\neq\emptyset$. So, $$k\leq\kappa'(G)\leq d_{G}(\bar{H}_2, V(G)\backslash \bar{H}_2)\leq |U|\leq k-1,$$ a contradiction.

If $H_1=\emptyset$ holds, then $\bar{H}_1\neq\emptyset$ by $F\neq\emptyset$. Thus $$k\leq\kappa'(G)\leq d_{G}(\bar{H}_1,V(G)\backslash\bar{H}_1)\leq |U|\leq k-1,$$  also a contradiction. $\Box$

\section{Concluding remarks}
\begin{thm}(Diwan and Tholiya \cite{Diwan}) For any tree $T$ with order $m$, every 1-connected graph $G$ with $\delta(G)\geq m$ contains a tree $T'\cong T$ such that $G-V(T')$ is 1-connected.
\end{thm}

Since the definitions of 1-edge-connected graph and 1-connected graph are equivalent, there exists a redundant subtree in a $1$-edge-connected graph with minimum degree of at least $m$ by Theorem 4.1.

\begin{thm} For any tree $T$ with order $m$, every 1-edge-connected graph $G$ with $\delta(G)\geq m$ contains a tree $T'\cong T$ such that $G-V(T')$ is 1-edge-connected.
\end{thm}

\begin{thm}(Hong et al. \cite{Hong}) For any tree $T$ with order $m$, every $2$-connected (respectively, $3$-connected) graph $G$ with $\delta(G)\geq m+2$ (respectively, $\delta(G)\geq m+3$) contains a subtree $T'\cong T$ such that $\kappa(G-V(T'))\geq 2$ (respectively, $\kappa(G-V(T'))\geq 3$).
\end{thm}

Using almost the same proof as Theorem 4.3, we can prove that there is a redundant subtree in a $2$-edge-connected (respectively, $3$-edge-connected) graph with minimum degree of at least $m+2$ (respectively, $m+3$). Specifically, we only need to substitute ``2-connected" and ``3-connected" in the proof of Theorem 4.3
with ``2-edge-connected" and ``3-edge-connected", respectively.

\begin{thm} For any tree $T$ with order $m$, every $2$-edge-connected (respectively, $3$-edge-connected) graph $G$ with $\delta(G)\geq m+2$ (respectively, $\delta(G)\geq m+3$) contains a subtree $T'\cong T$ such that $\kappa' (G-V(T'))\geq 2$ (respectively, $\kappa'(G-V(T'))\geq 3$).
\end{thm}

Motivated by the results above, we would like to propose the following conjecture, which relates to redundant subtree in a generally $k$-edge-connected graphs. As the case $k=1$ solved completely by Theorem 4.2, we assume $k\geq2$ in the following conjecture.

\begin{conj} For any tree $T$ of order $m$  and  any integer  $k\geq 2$, every $k$-edge-connected graph $G$ with minimum degree at least $k+m$ contains a subtree $T'\cong T$ such that $G-V(T')$ is still $k$-edge-connected.
\end{conj}

Since $K_{k+m}-V(T)$ is not $k$-edge-connected for any tree $T$ with order $m$, where $K_{k+m}$ denotes the complete graph of order $k+m$, the minimum degree bound would be best possible when Conjecture 2 were true.
Theorem 1.2 and Theorem 3.1 confirm Conjecture 2 for $m=1$ and $m=2$, respectively. For $k=2$ and $k=3$, by Theorem 4.4,   Conjecture 2 is also true.

{\flushleft\textbf{Declaration of competing interest}}

The authors declare that they have no known competing financial interests or personal relationships that could have appeared to influence the work reported in this paper.

%{\flushleft\textbf{Acknowledgements}}
%
%The authors would like to thank the editor and the anonymous reviewers  for their valuable and kind suggestions which greatly improved the original manuscript.

\end{document}